\documentclass[a4paper,12pt]{amsart}
\usepackage{amssymb}
\usepackage{ifthen}
\usepackage{graphicx}
\usepackage{float}
\usepackage{caption}
\usepackage{subcaption}
\usepackage{booktabs}
\usepackage[usenames]{color}
\usepackage{amssymb}
\usepackage{ifthen}
\usepackage{graphicx}
\usepackage{float}
\usepackage{caption}
\usepackage{subcaption}
\usepackage{booktabs}
\usepackage[usenames]{color}

\usepackage{latexsym,enumerate,amssymb, xcolor,mathrsfs}
\usepackage[colorlinks]{hyperref}

\newtheorem{thm}{Theorem}
\newtheorem{cor}{Corollary}
\newtheorem{lem}{Lemma}

\newtheorem{rem}{Remark}
\newtheorem{example}{Example}
\newtheorem{defn}{Definition}
\newtheorem{prob}{Problem}

\newtheorem{Conj}{Conjecture}

\newtheorem{Thm}{Theorem}

\newtheorem{Lem}{Lemma}

\theoremstyle{definition}

\newcounter {own}
\def\theown {\thesection       .\arabic{own}}

\newenvironment{pf}[1][]{%
	\vskip 3mm
	\noindent
	\ifthenelse{\equal{#1}{}}%
	{{\slshape Proof. }}%
	{{\slshape #1.} }%
}%
{\qed\bigskip}

\newcommand{\ID}{{\mathbb D}}

\def\be{\begin{equation}}
	\def\ee{\end{equation}}

\newcommand{\bee}{\begin{enumerate}}
	\newcommand{\eee}{\end{enumerate}}

\newcommand{\blem}{\begin{lem}}
	\newcommand{\elem}{\end{lem}}
\newcommand{\bthm}{\begin{thm}}
	\newcommand{\ethm}{\end{thm}}
\newcommand{\bcor}{\begin{cor}}
	\newcommand{\ecor}{\end{cor}}
\newcommand{\beg}{\begin{example}}
	\newcommand{\eeg}{\end{example}}
\newcommand{\begs}{\begin{examples}}
	\newcommand{\eegs}{\end{examples}}
\newcommand{\bdefe}{\begin{defn}}
	\newcommand{\edefe}{\end{defn}}
\newcommand{\bprob}{\begin{prob}}
	\newcommand{\eprob}{\end{prob}}
\newcommand{\bei}{\begin{itemize}}
	\newcommand{\eei}{\end{itemize}}

\newcommand{\bcon}{\begin{conj}}
	\newcommand{\econ}{\end{conj}}
\newcommand{\bcons}{\begin{conjs}}
	\newcommand{\econs}{\end{conjs}}
\newcommand{\bprop}{\begin{propo}}
	\newcommand{\eprop}{\end{propo}}
\newcommand{\br}{\begin{rem}}
	\newcommand{\er}{\end{rem}}
\newcommand{\brs}{\begin{rems}}
	\newcommand{\ers}{\end{rems}}
\newcommand{\bo}{\begin{obser}}
	\newcommand{\eo}{\end{obser}}
\newcommand{\bos}{\begin{obsers}}
	\newcommand{\eos}{\end{obsers}}
\newcommand{\bpf}{\begin{pf}}
	\newcommand{\epf}{\end{pf}}
\newcommand{\ba}{\begin{array}}
	\newcommand{\ea}{\end{array}}
\newcommand{\beq}{\begin{eqnarray}}
	\newcommand{\beqq}{\begin{eqnarray*}}
		\newcommand{\eeq}{\end{eqnarray}}
	\newcommand{\eeqq}{\end{eqnarray*}}

\newcounter{minutes}
\divide\time by 60
\newcounter{hours}
\multiply\time by 60 \addtocounter{minutes}{-\time}

\begin{document}
	\bibliographystyle{amsplain}
	\title[$K$-quasiconformal harmonic mappings]{Coefficients and Integral Mean Estimates for $K$-Quasiconformal Harmonic Mappings}
	
	\thanks{
		File:~\jobname .tex,
		printed: \number\year-\number\month-\number\day,
		\thehours.\ifnum\theminutes<10{0}\fi\theminutes}
	
	\author{Jasbir Parashar 
	}
	
	\address{Jasbir Parashar,
		Department of Mathematics,
		Indian Institute of Technology Ropar, Punjab--140001, India.}
	\email{jasbir.20maz0013@iitrpr.ac.in}
	
	\author{Saminathan Ponnusamy}
	\address{Saminathan Ponnusamy,
		Department of Mathematics,
		Indian Institute of Technology Madras, Chennai--600036, India.}
	\email{samy@iitm.ac.in}
	
	\author{A. Sairam Kaliraj}
	\address{A. Sairam Kaliraj, Department of Mathematics,
		Indian Institute of Technology Ropar, Punjab--140001, India.}
	\email{sairamkaliraj@gmail.com}
	
	
	\begin{abstract}
		Recently, Li and Ponnusamy~\cite{LiPonnusamy2025} established the coefficient conjecture proposed by Wang et al.~\cite{Wang2024} for several prominent geometric subclasses of $\mathcal{S}^0_H(K)$, the class of sense-preserving $K$-quasiconformal univalent harmonic mappings in the unit disk. In this paper, we show that the conjecture continues to hold for a class of $K$-quasiconformal harmonic mappings defined via quasi-subordination. Further, we determine the range of $p>0$ for which such mappings belong to the Hardy space ${\bf h}^p$ and the weighted Bergman space $\mathbf{a}^{\mathbf{p}}_{\boldsymbol{\beta}}$, for $\beta>-1$. Our Hardy space result makes significant progress toward a problem posed by Pavlovi\'{c}, while the Bergman space result sharpens the range obtained by Das and Rasila~\cite{DasRasila}, doubling the previously known bounds. In addition, we obtain refined growth and integral mean estimates for the subclass, improving earlier results and providing further evidence toward an open problem raised by Das et al.~\cite{DasRasila2025}. Parallel results are also discussed for odd $K$-quasiconformal harmonic mappings.
	\end{abstract}
	\thanks{ }

	\maketitle
	\pagestyle{myheadings}
	\markboth{J. Parashar, S. Ponnusamy and A. Sairam Kaliraj}{$K$-quasiconformal harmonic mappings}
	
	\section{Introduction and Problems}\label{pap4-sec1}
	\subsection*{1.1. Harmonic mappings.}
	We consider the class of complex-valued functions $f=u+iv$ on the unit disk $\mathbb{D} := \{ z \in \mathbb{C} : \vert z \vert < 1 \}$ that are \textit{harmonic}, i.e., $f$ satisfy $\Delta f = 0$, where $\Delta$
	denotes the Laplacian operator given by
	$$\Delta= 4 \partial_z \partial_{\bar{z}} = \partial_{xx} + \partial_{yy} \quad (z=x+iy).
	$$ Every harmonic function admits a decomposition of the form
	$f = h + \overline{g}$, where $h$ and $g$ are analytic in $\mathbb{D}$ and are respectively called analytic and co-analytic parts of $f$.
	It is convenient to work with functions that are normalized by the assumption $h(0)=g(0)=0$ and $h'(0)=1$, and thus have the  representation of the form
	\begin{equation}\label{equation1***}
		f(z)= z + \sum_{n=2}^{\infty} a_n z^n
		+ \overline{\sum_{n=1}^{\infty} b_n z^n},~\quad z \in \mathbb{D}.
	\end{equation}
	The class of all harmonic functions $f$ in $\mathbb{D}$ of the form \eqref{equation1***} is denoted by $\mathcal{H}$.
	We say that $f$ is \textit{sense-preserving} in $\mathbb{D}$ if its Jacobian $J_f(z)$ is positive for all $z \in \mathbb{D}$, where
	$$J_f(z) = \vert h'(z)\vert^2 - \vert g'(z)\vert^2 = \vert f_z(z)\vert^2 - \vert f_{\bar z}(z)\vert^2.
	$$
	Thus, every sense-preserving harmonic function $f$ induces an analytic function $w:\,\ID\to \ID$, known as the analytic  \textit{dilatation} of $f$, such that $g'(z) = w(z) h'(z)$.
	The subclass of $\mathcal{H}$ consisting of sense-preserving, univalent harmonic functions is denoted by $\mathcal{S}_H$. Although $\mathcal{S}_H$ is not compact, its subclass $\mathcal{S}_H^0 := \{ f \in \mathcal{S}_H : g'(0) = 0 \}$ is compact. However, both $\mathcal{S}_H$ and $\mathcal{S}_H^0$ are normal families.
	\subsection*{1.2.  $K$-quasiconformal Harmonic mappings.}
	A harmonic mapping $f = h + \overline{g}$ defined on $\mathbb{D}$ is said to be \emph{$K$-quasiregular} if it is sense-preserving and its \emph{analytic dilatation}
	$\omega = g'/h'$ satisfies the inequality
	$
	\vert \omega(z)\vert \le k < 1 ,
	$
	where
	$
	K= \frac{1+k}{1-k}\geqslant 1.
	$ Throughout the paper, we use $k$ and $K$ interchangeably, in view of this relation.
	We say that $f$ is a \emph{$K$-quasiconformal mapping} if it is a $K$-quasiregular and univalent in $\mathbb{D}$. Let us define
	$$
	\mathcal{S}_H(K) := \{ f  \in \mathcal{S}_H
	:\, f \text{ is } K\text{-quasiconformal} \}.
	$$
	We say that a function $f$ is \emph{harmonic quasiconformal} if $f\in \mathcal{S}_{H}(K)$ for some $K \ge 1$. Also, we define $$
	\mathcal{S}_H^0(K) := \left\{ f \in \mathcal{S}_H^0 : f \text{ is } K\text{-quasiconformal} \right\}.
	$$

	\subsection*{1.3. Hardy Spaces }
	For $0<p\le\infty$, the Hardy space $\bf{H}^p$ (resp.\ $\bf{h}^p$) consists of all complex-valued analytic (resp.\ harmonic) functions $\Phi$ on the unit disk $\mathbb{D}$ such that $\|\Phi\|_p<\infty$, where
	$$
	\|\Phi\|_p =
	\begin{cases}
		\sup\limits_{0<r<1} M_p(r,\Phi) & \text{if } p \in (0, \infty), \\
		\sup\limits_{z \in \mathbb{D}} \vert \Phi(z) \vert & \text{if } p = \infty,
	\end{cases}
	$$
	and the integral means $M_p(r,\Phi)$ are given by:
	$$
	M_p(r, \Phi) = \left( \frac{1}{2\pi} \int_0^{2\pi} \vert \Phi(re^{i\theta})\vert^p \, d\theta \right)^{1/p}.
	$$
	
	The boundary behaviour of univalent harmonic mappings has been a subject of considerable interest in geometric function theory. A central question in this area is to determine the precise range of $p>0$ for which such mappings belong to the harmonic Hardy space $\bf{h}^p$. Let us define
	$
	\alpha := \sup_{f \in \mathcal{S}_H} |a_2|.
	$
	The quantity $\alpha$ plays a fundamental role in determining the boundary behaviour of functions in $\mathcal{S}_H$.
	The study of the boundary behaviour of functions $f \in \mathcal{S}_H$ was initiated by Abu-Muhanna and Lyzzaik \cite{AbuMuhanna:Lyzzailk}, who proved that $f \in \bf{h}^p$ for
	$
	p < 1/(2\alpha+2)^2.
	$
	Later, Nowak \cite{Nowak} improved this range to
	$
	p < 1/\alpha^2,
	$
	and obtained sharp results showing that $f \in \bf{h}^p$ for $p < 1/2$ (respectively, $p < 1/3$) whenever $f$ is a convex (respectively, close-to-convex) harmonic function. These findings led her to conjecture that if $f \in \mathcal{S}_H$, then $f \in \bf{h}^p$ for all $p < 1/\alpha$ and that this order is sharp. In \cite{Das_Sairam Kaliraj_2024}, this conjecture was verified under an additional assumption on the analytic part of $f$. However, for the entire class $\mathcal{S}_H$, the conjecture remains open.
	For fundamental results on integral means and Hardy spaces, we refer to the
	comprehensive treatments in~\cite{Duren-book:H^P,Pavlovi}.
	\subsection*{1.4. Pavlovi\'{c}'s Problem on Harmonic Hardy spaces}
	It is known that every $K$-quasiconformal mapping of the unit disk belongs to the quasiconformal Hardy space $\bf{H}^p$ for all $p<1/(2K)$,
	and that this exponent range is sharp (cf.\ Theorem~3.2 in~\cite{AstalaKoskela}).
	Since univalent harmonic $K$-quasiconformal mappings are not only quasiconformal, but are also harmonic, the constant
	$1/(2K)$ appearing in Theorem~3.2 of \cite{AstalaKoskela} could be improved further for the class of harmonic $K$-quasiconformal mappings.
	Motivated by these observations, Pavlovi\'{c} \cite[p.~315]{Pavlovic:FunctionClasses} posed the following problem.
	\begin{prob}\label{pb1}
		Find the sharp
		$
		p_0 = p_0(K)
		$
		such that every  $K$-quasiconformal harmonic mapping belongs to $\bf{h}^p$ for $p < p_0$.
	\end{prob}
	We consider this problem in Section \ref{Sect-Main} and this result (see Theorem \ref{thm1++}) provides significant progress toward answering Pavlovi\'{c}'s problem.\\

	Initial progress on the integral mean estimate of $K$-quasiconformal harmonic mappings can be found in \cite{PonnusamyQiaoWang}.
	The authors in \cite{PonnusamyQiaoWang} defined the class $\mathcal{B}_H(\lambda)$ as
	$$
	\mathcal{B}_H(\lambda)
	= \left\{ f = h + \overline{g} \in \mathcal{H} : \|T_f\| \le 2\lambda \right\}
	$$
	where
	$$
	\|T_f\|
	:= \sup_{z \in \mathbb{D},\ \theta \in [0,2\pi]}
	(1 - \vert z \vert^2)
	\left \vert
	\frac{h''(z) + e^{i\theta} g''(z)}
	{h'(z) + e^{i\theta} g'(z)}
	\right \vert
	$$ denotes the pre-Schwarzian norm of  $h+e^{i\theta}g$ for each $\theta \in [0, 2\pi]$ (see \cite{HernandezMartin}) and proved that if $\lambda > 1$, then $\mathcal{B}_H(\lambda)  \cap \mathcal{S}_{H}(K) \subset \bf{h}^p$
	for every $0 < p < 1/(\lambda - 1)$.
	We note that the conclusion of this result can be obtained under a weaker hypothesis.
	In fact, it is sufficient to assume boundedness of the pre-Schwarzian norm of the analytic
	part $h$. To make this statement precise, we define the class
	$$
	\mathcal{F}_H(\lambda)
	:= \left\{ f = h + \overline{g} \in  \mathcal{S}_{H}(K) : \|\mathcal{P}_h\| \le 2\lambda \right\},
	$$
	where
	$$
	\|\mathcal{P}_h\|
	:= \sup_{z \in \mathbb{D}}
	~(1 - \vert z \vert^2)
	\left \vert
	\frac{h''(z)}
	{h'(z)}
	\right \vert
	$$
	denotes the pre-Schwarzian norm of the analytic part $h$ of $f$. We raise
	
	\bprob\label{prob2}
	Determine range of $p>0$ such that every $f \in \mathcal{F}_H(\lambda)$ belongs to harmonic Hardy space $\bf{h}^p$.
	\eprob
	
	We answer this problem in Section \ref{Sect-Main} (see Theorem \ref{thm4+-}).

	\subsection*{1.5. Bergman Spaces}
	For $-1<\beta<\infty$ and $0<p<\infty$,  the weighted harmonic Bergman spaces $\mathbf{a}^{\mathbf{p}}_{\boldsymbol{\beta}}$ are defined by
	$$
	{\mathbf{a}^{\mathbf{p}}_{\boldsymbol{\beta}}}
	:=
	\left\{
	f \in \mathcal{H}
	:\;
	\|f\|^p_{\mathbf{a}^{\mathbf{p}}_{\boldsymbol{\beta}}}
	=
	\int_{\mathbb{D}} \vert f(z)\vert^p \, dA_\beta(z)
	< \infty
	\right\},
	$$
	where the weighted area measure is
	$$
	dA_\beta(z)
	=
	(\beta+1)\bigl(1-\vert z \vert^2\bigr)^\beta\, dA(z),
	$$
	and $dA$ denotes the normalized area measure on $\mathbb{D}$ given by
	$$
	dA(z)=\frac{1}{\pi}\,dx\,dy
	=\frac{1}{\pi}\,r\,dr\,d\theta ,
	\qquad z=x+iy=re^{i\theta}\in\mathbb{D}.
	$$
	For the case $\beta=0$, we simply denote harmonic Bergman spaces as  $\bf{a^p}$. However, it is clear that ${\bf{h}^p} \subset \bf{a^p}$. For more  discussion on Bergman spaces, we refer to \cite{ChenPonnusamyWang2013,DurenSchuster, HKZ, KalajMestrovic2011}. We now raise the following
	
	\bprob
	Determine the range of $p$ such that every $f \in \mathcal{S}_H(K)$ belongs to the Bergman space $\mathbf{a}^{\mathbf{p}}_{\boldsymbol{\beta}}$.
	\eprob
	
	We solve this problem in Theorem \ref{thm7}.

	
	\subsection*{1.6. Coefficient Conjecture}\label{1.6}
	In 1984, Clunie and Sheil-Small~\cite{Clunie-Small-84} formulated a conjecture regarding coefficient estimates for normalized univalent harmonic functions.
	\begin{Conj}[Clunie and Sheil-Small]\label{conj:A*}
		Let $ f = h + \overline{g} \in \mathcal{S}_{H}^{0} $ have the series representation as in \eqref{equation1***}.
		Then for all integers $n \geq 2$, the following sharp estimates hold:
		\beqq\label{equation2}
		~~\big \vert \vert a_n \vert - \vert b_n\vert \big \vert \leq n, \quad \vert a_n \vert\leqslant \frac{(n+1)(2n+1)}{6} \quad\text{and} \quad
		\vert b_n \vert \leqslant \frac{(n-1)(2n-1)}{6}.
		\eeqq
		These bounds are attained by the harmonic Koebe function $K(z)$ defined by \begin{equation}\label{equation3}
			K(z) = \frac{z - \frac{1}{2}z^2 + \frac{1}{6}z^3}{(1-z)^3} + \overline{\left( \frac{\frac{1}{2}z^2 + \frac{1}{6}z^3}{(1-z)^3} \right)}.
		\end{equation}
	\end{Conj}
	This conjecture is viewed as the harmonic analogue of the classical Bieberbach conjecture for the class
	$$
	\mathcal{S}
	=
	\{ f = h+\overline{g} \in \mathcal{S}_H : g \equiv 0 \text{ in } \mathbb{D} \}.
	$$
	For the class $\mathcal{S}$, de~Branges~\cite{L. de Branges} proved that $\vert a_n \vert=\vert f^{(n)}(0)/n! \vert\le n$ for all $n\ge2$ and for $f\in\mathcal{S}$, thereby settling the Bieberbach conjecture. A comprehensive discussion on the class $\mathcal{S}$ can be found from the books of  Duren \cite{Duren-book1} and Pommerenke \cite{Pommerenke}.
	
	
	The coefficient conjecture of Clunie and Sheil-Small has been verified for several geometric subclasses of $\mathcal{S}_H^{0}$, including the class $\mathcal{S}_H^{*0}$ of starlike functions,  the class $\mathcal{T}_H^{0}$  of typically real functions, the class $\mathcal{K}_H^{0}$ of close-to-convex functions and consequently for the class  $\mathcal{K}_H^{0}(\theta)~(0 \le \theta < \pi)$ of functions convex in the  direction of $\theta$ (see \cite{Clunie-Small-84,Duren:Harmonic,Precise}). Nevertheless, the conjecture remains unresolved for the entire class of univalent harmonic mappings. Much like the Bieberbach conjecture played a central role in the development of classical univalent function theory, the Clunie and Sheil-Small's coefficient conjecture has served an analogous role in advancing the theory of planar univalent harmonic mappings.
	Ponnusamy and Sairam Kaliraj~\cite{PonKal2015} made a significant contribution in this direction by establishing the coefficient estimates for the subclass:
	$$
	\mathcal{S}^0_H(\mathcal{S}) = \left\{ h + \overline{g} \in \mathcal{S}^0_H : h + e^{i\theta}g \in \mathcal{S} \text{ for some } \theta \in \mathbb{R} \right\}
	.$$
	Although it was initially conjectured that $\mathcal{S}^0_H(\mathcal{S}) = \mathcal{S}^0_H$, this equality was later shown to be false by Bshouty et al.~\cite{Bshouty}, highlighting the subtle nature of the relationship between these classes.
	A notable advancement in this direction was made in~\cite{Jasbir} using the framework of quasi-subordination.
	Recall that for analytic functions $\phi$ and $\psi$ in $\ID$, we say that the function $\phi$ is \emph{quasi-subordination} to the function $\psi$ (denoted $\phi \prec_q \psi$) if there exist analytic functions $\xi$ and $\eta$ with $\vert \xi(z)\vert \leq 1$, $\vert \eta(z)\vert < 1$, and $\eta(0) = 0$, such that:
	\[
	\phi(z) = \xi(z) \psi(\eta(z)), \quad z\in \ID.
	\]
	This notion generalizes both \emph{subordination} and \emph{majorization}: it reduces to subordination when $\xi\equiv1$, and to majorization when $\eta(z)=z$. For further details, see~\cite{Pommerenke, Robertson}.
	In ~\cite{Jasbir}, the authors established the validity of Conjecture~\ref{conj:A*} for the class $\mathcal{S}_{H_q}^{0}(\mathcal{S})$ defined as
	\begin{align*}
		\mathcal{S}_{H_q}^{0}(\mathcal{S}) &= \left\{
		f = h + \overline{g} \in \mathcal{S}_{H}^{0} :
		h + e^{i\theta}g \prec_q \psi \text{ for some } \theta \in \mathbb{R} \text{ and } \psi \in \mathcal{S}
		\right\}.
	\end{align*}
	This class contains all functions in $\mathcal{S}_{H}^{0}(\mathcal{S})$ and possibly more, as it relaxes the condition $h+e^{i\theta}g\in\mathcal{S}$ to a quasi-subordination relation.
	Next, we discuss a similar type of development for $K$-quasiconformal harmonic mappings.

	More recently, Wang \emph{et al.}~\cite{Wang2024} constructed a $K$-quasiconformal harmonic Koebe-type function $\mathbf{P}_{k}(z)$, defined by
	\begin{align}\label{EQ4*}
		\mathbf{P}_{k}(z)
		&= \mathbf{H}_{k}(z) + \overline{\mathbf{G}_{k}(z)} \notag \\[4pt]
		&= \frac{1}{(k-1)^3}
		\left(
		\frac{(k-1)(1-3k+2kz)z}{(1-z)^2}
		+ k(k+1)\log\!\left(\frac{1-z}{1-kz}\right)
		\right) \notag \\[6pt]
		&\quad
		+
		\frac{k}{(k-1)^3}\overline{
			\left(
			\frac{(1-k)(1+k-2z)z}{(1-z)^2}
			+ (k+1)\log\!\left(\frac{1-z}{1-kz}\right)
			\right)
		} \notag \\[6pt]
		&= z + \sum_{n=2}^\infty \mathbf{A}(n,k)\, z^n
		+ \sum_{n=2}^\infty \mathbf{B}(n,k)\, \overline{z}^{\,n},
	\end{align}
	where 
	\begin{equation}\label{+5}
		\mathbf{A}(n,k) = \frac{1}{n(1 - k)^{3}}
		\left( n^{2} + (-2n^{2} - 2n + 1)k + (n + 1)^{2}k^{2} - k^{n+1} - k^{n+2} \right)
	\end{equation}
	and
	\begin{equation}\label{+6}
		\mathbf{B}(n,k) = \frac{k}{n(1 - k)^{3}}
		\left( (n - 1)^{2} + (-2n^{2} + 2n + 1)k + n^{2}k^{2} - k^{n} - k^{n+1} \right),
	\end{equation}
	respectively.
	In connection with Problem~\ref{pb1} and motivated by the harmonic $K$-quasiconformal
	Koebe-type function $\boldsymbol{P}_{k}(z)$, the authors of \cite{Wang2024}
	proposed several related conjectures. One of these conjectures reads as follows.
	\begin{Conj}\label{conj:A}
		Suppose that $f = h + \overline{g} \in \mathcal{S}_{H}^{0}(K)$ is of the form  \eqref{equation1***}. Then
		\begin{equation}\label{EQ9*}
			\big\vert\,\vert a_{n}\vert - \vert b_{n}\vert\,\big\vert \leq n,
			\quad \vert a_{n}\vert \leq \mathbf{A}(n,k),
			\quad \text{and} \quad \vert b_{n}\vert \leq \mathbf{B}(n,k),
			\quad \text{for } n = 2,3,\ldots,
		\end{equation}
		where $\mathbf{A}(n,k)$ and $\mathbf{B}(n,k)$ are defined by \eqref{+5} and \eqref{+6}, respectively.
		Equality holds for the function $\boldsymbol{P}_{k}$ given by \eqref{EQ4*}.
	\end{Conj}
	Motivated by the progress made for the class $\mathcal{S}_{H}^{0}$, Li and Ponnusamy \cite{LiPonnusamy2025} very recently verified this conjecture for several geometric subclasses of $\mathcal{S}_H^{0}(K)$,
	including  $\mathcal{S}_H^{*0}(K)$, $\mathcal{K}_H^{0}(K)$,   $\mathcal{T}_H^{0}(K)$ and an analog coefficient inequalities for the class $\mathcal{C}_H^{0}(K)$ of $K$-quasiconformal harmonic convex functions,
	where
	$$
	\left \{\ba{rl}
	{\mathcal C}_H(K)&:={\mathcal S}_H(K)\cap {\mathcal C}_H ~\mbox{ and }~{\mathcal C}_H^0(K):={\mathcal S}_H(K)\cap {\mathcal C}_H^0;\\
	{\mathcal S}_H^{*}(K)&:={\mathcal S}_H(K)\cap {\mathcal S}_H^* ~\mbox{ and }~{\mathcal S}_H^{*0}(K):={\mathcal S}_H(K)\cap {\mathcal S}_H^{*0};\\
	{\mathcal K}_H(K)&:={\mathcal S}_H(K)\cap {\mathcal K}_H ~\mbox{ and }~{\mathcal K}_H^{0}(K):={\mathcal S}_H(K)\cap {\mathcal K}_H^0.
	\ea\right .
	$$
	Here ${\mathcal T}_H(K)$ is the class of typically real $K$-quasiconformal harmonic mappings and (cf. \cite{LiPonnusamy2025})
	$${\mathcal T}_H^0(K):={\mathcal T}_H(K)\cap {\mathcal T}_H^0.
	$$
	Let us now introduce
	$$
	\mathcal{S}^0_{H_q}(\mathcal{S},K)
	:=
	\left\{
	f=h+\overline{g} \in \mathcal{S}^0_H(K)
	:\;
	h+e^{i\theta}g \prec_q \psi
	\ \text{for some } \theta \in \mathbb{R}
	\ \text{and } \psi \in \mathcal{S}
	\right\}.
	$$
	This class may be regarded as a natural counterpart of the class $\mathcal{S}_{H_q}^{0}(\mathcal{S})$ in the setting of $\mathcal{S}^0_H(K)$. Denote by
	$\mathcal{S}_{H_s}^{0}(\mathcal{S},K)$ when the quasi-subordination symbol $\prec_q$ in the definition of $\mathcal{S}^0_{H_q}(\mathcal{S},K)$ is replaced just by subordination symbol $\prec$.
	It is worth pointing out that the class $\mathcal{S}_{H_s}^{0}(\mathcal{S},K)\cap  \mathcal{S}$ has been well-studied class introduced by Robertson \cite{Robertson,Robertson3} (see also
	Rogosinski \cite{Rogo-43}) and investigated by a number of authors.

	Motivated by the above discussion and the recent investigation on this topic, it is natural to raise the following.
	
	\bprob
	Does Conjecture \ref{conj:A} hold for $\mathcal{S}_{H_{q}}^0(\mathcal{S},K)$?
	\eprob
	
	In our next result (see Theorem \ref{thm1+}), we confirm  that Conjecture~\ref{conj:A} holds for the class $\mathcal{S}_{H_q}^{0}(\mathcal{S},K)$. In addition,  we obtain corresponding growth estimates.

	In~\cite{DasRasila2025}, Das et al. determined the range of $p>0$ for which
	$K$-quasiconformal convex and close-to-convex harmonic mappings belong to the
	Hardy space $\bf{h}^p$.
	More precisely, they proved that $K$-quasiconformal convex harmonic mappings
	belong to $\bf{h}^p$ for all $p<1$, while $K$-quasiconformal close-to-convex harmonic
	mappings belong to $\bf{h}^p$ for all $p<1/2$.
	These results refine the sharp bounds obtained by Nowak~\cite{Nowak}, who showed
	that convex harmonic mappings belong to $\bf{h}^p$ for all $p<1/2$, whereas
	close-to-convex harmonic mappings belong to $\bf{h}^p$ for all $p<1/3$.
	One may therefore observe that quasiconformal mappings belonging to the
	geometric subclasses of $\mathcal{S}_H(K)$ have restricted growth, as expected. Motivated by these observations, the authors of~\cite{DasRasila2025} posed the
	following question:
	
	\bprob\label{prob-hardy}
	Does every function $f\in\mathcal{S}_H(K)$ belong to $\bf{h}^p$
	for $p<1/2$?
	\eprob
	
	In our next result, we show that (see Theorem \ref{thm2@}) this question has an affirmative answer for the
	class $\mathcal{S}_{H_q}(\mathcal{S},K)$, where
	$$
	\mathcal{S}_{H_q}(\mathcal{S},K) := \left\{ h + \overline{g} \in \mathcal{S}_{H}(K) : h+e^{i\theta}g \prec_q \psi
	\ \text{for some } \theta \in \mathbb{R}
	\ \text{and } \psi \in \mathcal{S}\right\}.
	$$
	This extends and refines earlier work of Kayumov \emph{et al.}~\cite{Kayumov},
	who proved that if $f\in\mathcal{S}_{H}^{0}(\mathcal{S})$, then
	$f\in \bf{h}^p$ for all $p<1/3$.  In fact in Theorem \ref{thm2@}, we show that every $f \in \mathcal{S}_{H_q}(\mathcal{S}, K)$ belongs to $\mathbf{a}^{\mathbf{p}}_{\boldsymbol{\beta}}$ whenever $0<p<(2+\beta)/2$ and $\beta >-1$.

	At this place, an interesting comparison can be made with a classical result of Astala and Koskela~\cite{AstalaKoskela}.
	Theorem~3.2 in~\cite{AstalaKoskela} states that a $K$-quasiconformal mapping of the unit disk belongs to the quasiconformal Hardy space $\bf{H}^p$ for $p < 1/(2K)$, and that this range of $p$ is optimal.
	However, since $1/(2K) \ll 1/2$ when $K \gg 1$, the admissible range of $p$ obtained in Theorem~\ref{thm2@} is strictly larger.
	This observation indicates that investigating Hardy spaces associated with harmonic quasiconformal mappings may yield new and non-trivial results.

	\subsection*{1.7. Odd
		$K$-quasiconformal  harmonic mappings}\label{section1.7}
	The class of odd univalent analytic functions was also instrumental in  settling the celebrated Bieberbach conjecture~\cite{L. de Branges, Milin, MilinEstimate, Robertson2}.
	Motivated by the unresolved harmonic analogue of the Bieberbach conjecture, the authors in \cite{Jaglan} investigated several subclasses of odd univalent harmonic mappings in order to explore a possible analogous connection between odd univalent harmonic mappings and the harmonic Bieberbach conjecture. In view of Conjecture~\ref{conj:A}, we therefore focus on the study of odd
	$K$-quasiconformal  harmonic mappings.
	Our next result (see Theorem \ref{thm5}) concerns in establishing  coefficient and growth estimates for odd
	$K$-quasiconformal harmonic mappings in the class $\mathcal{S}^0_H(\mathcal{S},K)$, where
	$$
	\mathcal{S}^0_H(\mathcal{S},K) := \left\{ h + \overline{g} \in \mathcal{S}^0_H(K) : h + e^{i\theta}g \in \mathcal{S} \text{ for some } \theta \in \mathbb{R} \right\}.
	$$

	In our final result (see Theorem \ref{thm8++}) we determine the range of $p$ for which odd functions in the class
	$S_{H}(\mathcal{S},K)$ belong to the Hardy space $\bf{h}^{p}$ and the Bergman space $\mathbf{a}^{\mathbf{p}}_{\boldsymbol{\beta}}$, where
	$$
	\mathcal{S}_H(\mathcal{S},K) := \left\{ h + \overline{g} \in \mathcal{S}_H(K) : h + e^{i\theta}g \in \mathcal{S} \text{ for some } \theta \in \mathbb{R} \right\}.
	$$
	This result refines an earlier theorem of Jaglan and Sairam \cite{Jaglan}, which asserts that if $f$ is an odd function in the class
	$S_{H}^{0}(\mathcal{S})$, then $f \in \bf{h}^{p}$ for $0<p<1/2$.


	\section{Preliminaries and  Auxiliary Results}\label{pap2-sec1.2}
	This section presents foundational results from the literature that are required to prove our main theorems, beginning with a generalization of the de Branges' theorem for the class $\mathcal{S}$ within the framework of quasi-subordination.
	
	\begin{Thm}\cite[p.249]{Hayman}\label{thmA}
		Let $\xi(z)$ and $\eta(z)$ be analytic functions in the unit disk $\mathbb{D}$ satisfying $\vert \xi(z) \vert \leq  1$, $\vert \eta(z) \vert < 1$ for all $z \in \mathbb{D}$, with $\eta(0) = 0$. If $\psi \in \mathcal{S}$ and
		$$
		\phi(z) = \xi(z)\psi(\eta(z)) = \sum_{n=1}^\infty A_n z^n,
		$$
		then the coefficients satisfy $\vert A_n \vert \leq n$ for all $n \geq 1$. Equality holds if and only if $\phi$ is a rotation of the Koebe function $
		\boldsymbol{k}(z) = z/(1 - z)^2$.
	\end{Thm}
	The following result of Clunie and Sheil-Small~\cite{Clunie-Small-84}, known as the shearing method, provides a criterion for the univalency of harmonic mappings and serves as a fundamental tool for constructing extremal functions.
	
	\begin{Thm}\label{clunie}
		Let $f = h + \overline{g}$ be harmonic and locally univalent in the unit disk $\mathbb{D}$.
		Then $f$ is univalent and its range is convex in the direction of $\alpha$
		if and only if $h - e^{2 i \alpha} g$ has the same properties.
	\end{Thm}
	For each  $f \in \mathcal{S}$, the odd function $\phi$ defined by (the square-root transform of $f$)
	$$
	\phi(z)= z\left(f(z^2)/z^2\right)^{1/2}
	= z + \sum_{n=1}^{\infty} a_n z^{2n+1}
	$$
	belongs to $\mathcal{S}$, and conversely, every odd univalent analytic function admits a representation as a square-root transformation of some function $f\in\mathcal{S}$. Let $\mathcal{S}^{2}$ denote the class of all odd functions in the class $\mathcal{S}$.
	The following result establishes distortion theorem for the class $\mathcal{S}^{2}$, which is crucial in the proof of Theorem \ref{thm5}.
	
	\begin{Lem}\cite{Jasbir}\label{Lemma5A}
		{\it
			Let $\phi \in \mathcal{S}^{2}$. Then the following inequalities hold:
			\begin{align*}
				\frac{1 - r^2}{\left(1 + r^2\right)^{2}}
				&\leq \lvert \phi'(z) \rvert
				\leq \frac{1 + r^2}{\left(1 - r^2\right)^{2}}~\quad \text{for} ~\quad \vert z \vert = r < 1.
			\end{align*}
			In addition, sharpness of the upper and lower bounds is achieved by $z/(1-z^2)$ and $z/(1+z^2)$, respectively.
		}
	\end{Lem}
	
	The following result provides a universal bound for the coefficients of odd univalent functions. 
	
	\begin{Thm}\label{thmA1}\cite{Hu}
		Let $f \in \mathcal{S}^{2}$ and  $f(z) = z + \sum_{n=1}^{\infty} a_{n} z^{2n+1}$. Then $\vert a_{n} \vert < 1.1305$ for all $n \geq 1$.
	\end{Thm}
	
	
	The following lemma provides key technical tools for establishing integral mean estimates.
	
	\begin{Lem}\cite{PonnusamyQiaoWang}\label{lemB++}
		{\it
			Let $f \in \mathcal{S}_{H}(K) $ and $p > 0$. Then
			$$
			M_p^{p}(r,f)
			\le
			\frac{2(1+k^2)\bigl(\vert p-2 \vert + 1\bigr)}{1-k^2}
			\int_0^r M(\rho)^p \rho^{-1}\, d\rho,
			\qquad 0 \le r < 1,
			$$
			where
			$$
			M(r) := M(r,f)
			=
			\max_{0 \le \theta \le 2\pi}
			\vert f(re^{i\theta}) \vert .
			$$
		}
	\end{Lem}
	\
	\section{Main Results and their proofs}\label{Sect-Main}
	
	As an answer to Problem \ref{pb1}, we state and prove our first main result which determines the range of $p=p(K)$ such that every $K$-quasiconformal harmonic mapping belongs to the Hardy space $\bf{h}^{p}$.
	
	\begin{thm}\label{thm1++}
		Let $f = h + \overline{g} \in \mathcal{S}_H(K)$. Then $f\in \bf{h}^p$ for
		$0 < p < 1/\alpha_K$, where
		\beq \label{eqUation}
		\alpha_K := \sup_{f \in \mathcal{S}_H(K)} \left \vert \frac{h''(0)}{2} \right \vert.
		\eeq
	\end{thm}
	\bpf 
	Let $f = h + \overline{g}$ belong to the class $\mathcal{S}_H(K)$ and have the form \eqref{equation1***}. Fix $\zeta \in \mathbb{D}$ and apply a disk automorphism to obtain the function
	$$
	F(z)
	=
	\frac{
		f\!\left( \dfrac{z+\zeta}{1+\overline{\zeta}z} \right) - f(\zeta)
	}{
		(1-\vert \zeta \vert^2) h'(\zeta)
	}= H(z)+\overline{G(z)},
	\qquad z \in \mathbb{D}.
	$$
	Note that the class $\mathcal{S}_H(K)$ is invariant under disk automorphisms (cf. \cite[p.6]{Duren:Harmonic}),  that is, $F \in \mathcal{S}_H(K)$. A routine calculation shows that
	$$
	H''(0)
	=
	(1-\vert \zeta \vert^2)\frac{h''(\zeta)}{h'(\zeta)} - 2\overline{\zeta}.
	$$
	Since $\vert H''(0) \vert \le 2\alpha_K$, it follows that
	$$
	\operatorname{Re}\!\left\{\frac{z h''(z)}{h'(z)}\right\}
	\le
	\frac{2r^2 + 2\alpha_K r}{1-r^2},
	\qquad \vert z \vert = r < 1.
	$$
	The last inequality can be rewritten in the form
	$$
	\frac{\partial}{\partial r}
	\bigl(\log \vert h'(r e^{i\theta}) \vert\bigr)
	\le
	\frac{2r + 2\alpha_K}{1-r^2}.
	$$
	Now, a methodology similar to that of   \cite[p.~98]{Duren:Harmonic} leads to the following estimate
	$$
	\vert h'(z) \vert
	\le
	\frac{(1+r)^{\alpha_K-1}}{(1-r)^{\alpha_K+1}},
	\qquad \vert z \vert = r < 1.
	$$
	Since $\vert \omega(z) \vert \leq k$, where $k= \frac{K-1}{K+1}$, we obtain
	$$
	\vert g'(z)\le k \vert h'(z)  \vert
	\le
	\frac{K-1}{K+1}\cdot\frac{(1+r)^{\alpha_K-1}}{(1-r)^{\alpha_K+1}},
	\qquad \vert z \vert = r < 1,
	$$
	so that
	$$ \vert h'(z) \vert + \vert g'(z) \vert \le \frac{2K}{K+1}\cdot\frac{(1+r)^{\alpha_K-1}}{(1-r)^{\alpha_K+1}}.$$
	Next, we have
	\begin{align}
		\vert f(z) \vert
		&=
		\left\vert
		\int_{\Gamma}
		\frac{\partial f}{\partial \zeta}\, d\zeta
		+
		\frac{\partial f}{\partial \overline{\zeta}}\, d\overline{\zeta}
		\right\vert \notag
		\le
		\int_{\Gamma}
		\bigl( \vert h'(\zeta) \vert + \vert g'(\zeta) \vert \bigr)\, \vert d\zeta \vert \notag \\
		&\le
		\frac{2K}{K+1}
		\int_0^r
		\frac{(1+s)^{\alpha_K-1}}{(1-s)^{\alpha_K+1}}\, ds \notag
		=
		\frac{K}{K+1} \cdot \frac{1}{\alpha_K}
		\left[
		\left( \frac{1+r}{1-r} \right)^{\alpha_K}
		- 1
		\right],
	\end{align}
	for $\vert z \vert = r < 1$, where $\Gamma$ is the radial line segment from $0$ to $z$. Now, invoking Lemma~\ref{lemB++}, we obtain
	\begin{align*}
		\lim_{r \to 1^-} M_p^p(r,f)
		&\le \frac{2(1+k^2)(\lvert p-2\rvert+1)}{1-k^2}
		\int_0^1 M(s)^p s^{-1}\, ds \\
		&\le \frac{(1+K^2)(\lvert p-2\rvert+1)}{K}
		\left( \frac{K}{\alpha_K (K+1)} \right)^p \notag \\
		&\quad \times \int_0^1
		\frac{\big((1+s)^{\alpha_K} - (1-s)^{\alpha_K}\big)^p}
		{s(1-s)^{\alpha_K p}} \, ds \\
		&\le C(p,K) \int_0^1
		\frac{s^{p-1}}{(1-s)^{\alpha_K p}} \, ds = C(p,K)\, B(p, 1-\alpha_K p),
	\end{align*}
	where $B(p,1-\alpha_Kp)$ denotes the beta function and  $C(p,K)> 0$ is a constant depending on $p$ and $K$. Therefore, $f \in {\bf{h}^p} \text{ for } 0<p < 1/\alpha_K.$
	\epf
	
	\begin{rem}
		By comparing Theorem~\ref{thm1++}  with Theorem~3.2 in~\cite{AstalaKoskela} and Conjecture~2 in~\cite{Wang2024}, we observe that the range $0 < p < 1/(2K)$ is not sharp for $K$-quasiconformal harmonic mappings when $K$ is large. Indeed, $ 1/(2K)$ approaches $0$ for large $K$, whereas $1/\alpha_K \ge 1/\alpha$ (since $\alpha_K\le \alpha$). Consequently, Theorem~\ref{thm1++} provides a better estimate for sufficiently large values of $K$. Moreover, since the class $\mathcal{S}\subset \mathcal{S}_H(K)$ for $K\ge 1$, it follows that $\alpha_K \ge 2$. Therefore, we have
		$$
		\frac{1}{\alpha} \le \frac{1}{\alpha_K} \le \frac{1}{2}.
		$$
		If sharpness can be established in Theorem~1, or if one can prove that the range $0<p<1/2$ holds in Theorem~1, then Pavlovi\'{c}'s problem would be completely answered.
	\end{rem}

	We state and prove our next result which determines the range of $p>0$ such that every $f \in \mathcal{F}_H(\lambda)$ belongs to harmonic Hardy space $\bf{h}^p$. This answers Problem \ref{prob2}.
	
	\begin{thm}\label{thm4+-}
		For $\lambda >1$, let $f = h + \overline{g} \in \mathcal{F}_H(\lambda)$.  Then $f\in \bf{h}^p$ for $0 < p < 1/(\lambda - 1)$.
	\end{thm}

	\bpf
	Let $f = h + \overline{g}\in \mathcal{F}_H(\lambda)$ and have the form \eqref{equation1***}. Then $h$ is analytic in $\mathbb{D}$ with $h(0)=0$, $h'(0)=1$, and
	\beqq
	(1-\vert z \vert^{2})\left \vert \frac{h''(z)}{h'(z)}\right \vert \le 2\lambda, \qquad z\in\mathbb{D}.
	\eeqq
	Applying an argument similar to that used in \cite[Theorem~3.1]{BeckerPommerenke}, we obtain
	\beqq
	\left(\frac{1-\vert z \vert}{1+\vert z \vert}\right)^{\!\lambda}
	\le \vert h'(z)\vert \le
	\left(\frac{1+\vert z \vert}{1-\vert z \vert}\right)^{\!\lambda}.
	\eeqq
	Since  $\vert \omega(z) \vert \leq k$, as in the proof of Theorem~1, it follows that
	\begin{align}
		\vert f(z) \vert
		&
		\le
		\int_{\Gamma}
		\bigl( \vert h'(\zeta) \vert + \vert g'(\zeta) \vert \bigr)\, \vert d\zeta \vert \notag \\
		&\le
		\frac{2K}{K+1}
		\int_0^r
		\frac{(1+s)^{\lambda}}{(1-s)^{\lambda}}\, ds \notag
		\le
		\frac{2K}{K+1} \cdot \frac{2^{\lambda}}{\lambda-1} \cdot
		\frac{1}{(1-r)^{\lambda-1}},
	\end{align}
	for $\vert z \vert = r < 1$, where $\Gamma$ is the radial line segment from $0$ to $z$. Now, invoking Lemma~\ref{lemB++}, we obtain the desired result.
	\epf
	
	Note that the function $Q_k$, defined by (cf. \cite[Eq.~(1.10)]{LiPonnusamy2025})
	$$
	Q_k(z) = \boldsymbol{P}_{k}(z) + k\,\overline{\boldsymbol{P}_{k}(z)},
	$$
	belongs to the class $\mathcal{S}_{H}(K)$, where $\boldsymbol{P}_{k}(z)$ is defined  in \eqref{EQ4*}. This function serves as an extremal function for several subclasses of the class $\mathcal{S}_{H}(K)$ (see \cite[Corollary~1.1]{LiPonnusamy2025}). In the limiting case $k \to 1^{-}$, the function $Q_k(z)$ coincides with
	$$
	Q(z) = K(z) +\overline{ b_1 K(z)}, \qquad \vert b_1 \vert < 1,
	$$
	where $K(z)$ is harmonic Koebe function defined by \eqref{equation3}.
	It is known that $Q \in \bf{h}^{p}$ for $0<p<1/3$ (see \cite{Nowak}).
	A similar conclusion follows from Theorem~\ref{thm4+-} and Table~1.
	The values of $\|\mathcal{P}_{\mathcal{H}_k}\|$ (correct up to two decimal places) for various values of $k$ are presented in Table~1, where $\mathcal{H}_k$ denotes the analytic part of the function $Q_k$ defined above. Moreover, it follows that for $0 < k \leq 0.45$, the function $Q_k$ belongs to $\bf{h}^{p}$ for all $0<p<1/2$, thereby providing supporting evidence for the Problem \ref{prob-hardy}.

	\begin{table}[h!]
		\centering
		\renewcommand{\arraystretch}{1.2} 
		\begin{tabular}{c @{\hspace{1.5cm}} c @{\hspace{1.5cm}} c}
			\toprule
			$k$ & $\|\mathcal{P}_{\mathcal{H}_{k}}\|$ & $\lambda$ \\
			\midrule
			0.10 & 6.00 & 3.00 \\
			0.20 & 6.00 & 3.00 \\
			0.30 & 6.00 & 3.00 \\
			0.40 & 6.00 & 3.00 \\
			0.50 & 6.01 & 3.00 \\
			0.60 & 6.12 & 3.06 \\
			0.70 & 6.30 & 3.15 \\
			0.80 & 6.56 & 3.28 \\
			0.90 & 6.94 & 3.47 \\
			0.9999 & 7.96 & 3.98\\
			\bottomrule
		\end{tabular}
		\caption{Approximate values of $\|\mathcal{P}_{\mathcal{H}_{k}}\|$ for different values of $k$}
	\end{table}
	
	We next determine the range of $p$ such that every $f \in \mathcal{S}_H(K)$ belongs to the Bergman space $\mathbf{a}^{\mathbf{p}}_{\boldsymbol{\beta}}$. More precisely, we prove the following:
	
	\begin{thm}\label{thm7}
		We have the inclusion $\mathcal{S}_H(K)\subset    \mathbf{a}^{\mathbf{p}}_{\boldsymbol{\beta}}$ for
		$0 < p < (\beta+2)/\alpha_K$, where $\alpha_K$ is defined in \eqref{eqUation}.
	\end{thm}
	\bpf 
	Let $f \in \mathcal{S}_H(K)$ and $p>0$. Then, we have
	\beqq
	\vert\vert f \vert\vert_{\mathbf{a}^{\mathbf{p}}_{\boldsymbol{\beta}}}^p
	=\int_{\mathbb{D}} \vert f(z)\vert^p \, dA_\beta(z)&=&(\beta+1)\int_{\mathbb{D}} \vert f(z)\vert^p (1- \vert z \vert^2)^{\beta}\, dA(z)\\
	&=&
	2(\beta+1)\int_0^1 r(1-r^2)^{\beta} M_p^p(r,f)\, dr
	\eeqq
	Now applying Lemma~\ref{lemB++}, we obtain
	\beqq
	\vert\vert f \vert\vert_{\mathbf{a}^{\mathbf{p}}_{\boldsymbol{\beta}}}^p
	&\le &
	\frac{4(\beta+1)(1+k^2)\bigl(\vert p-2 \vert + 1\bigr)}{1-k^2} \int_0^1 r(1-r^2)^{\beta}\left(\int_0^r M(\rho)^p \rho^{-1}\,d\rho\right) dr \\
	&=&
	\frac{4(\beta+1)(1+k^2)\bigl(\vert p-2 \vert + 1\bigr)}{1-k^2}\int_0^1 M(\rho)^p \rho^{-1}
	\left(\int_\rho^1 r(1-r^2)^{\beta}\,dr\right) d\rho \\
	&\lesssim&
	\int_0^1 M(\rho)^p \rho^{-1}(1-\rho^2)^{\beta+1}\,d\rho\\
	&\lesssim&\int_{0}^{1}
	\frac{\left((1+\rho)^{\alpha_K}-(1-\rho)^{\alpha_K}\right)^p}
	{\rho(1-\rho)^{\alpha_K p}}(1-\rho^2)^{\beta+1}
	\, d\rho\\
	&\lesssim& \int_0^1 (1-\rho)^{-\alpha_K p+\beta+1}\, \rho^{p-1}\,d\rho.
	\eeqq
	The last integral is the beta function $B(p, 2+\beta-\alpha_K p)$ and converges for all $0<p<(\beta+2)/\alpha_K$.
	Therefore,
	$
	f \in \mathbf{a}^{\mathbf{p}}_{\boldsymbol{\beta}}$
	~whenever~ $
	0 < p < (\beta+2)/\alpha_K.
	$
	\epf
	
	\begin{rem} In the special case $\beta = 0$,  Theorem \ref{thm7} refines the range obtained by Das and Rasila~\cite[Theorem~3]{DasRasila}, yielding precisely twice the earlier range.
	\end{rem}
	
	Next we state and prove the following result which contains solution to Conjecture~\ref{conj:A}  for the class $\mathcal{S}_{H_q}^{0}(\mathcal{S},K)$.
	
	\begin{thm}\label{thm1+}
		Let $f = h + \overline{g} \in \mathcal{S}_{H_{q}}^0(\mathcal{S},K)$ and have the series representation of the form \eqref{equation1***}. Then for all integers $n \geq 2$, the coefficients of $h$ and $g$ satisfy the sharp inequalities in \eqref{EQ9*}. Furthermore,
		for \( \vert z \vert = r < 1 \),
		\beqq
		\vert f(z)\vert \leq S(k,r),
		\eeqq
		where
		$$S(k,r)=\frac{r (k + 1)}{( 1-k)( 1-r)^2}- \frac{4 k r}{(1-k)^2 (1-r)}
		+ \frac{2 k (1 + k)}{(k - 1)^3}
		\log\left( \frac{ 1-r}{1-k r} \right) .
		$$
		These estimates are sharp for the class $\mathcal{S}_{H_{q}}^0(\mathcal{S},K)$, and equality is attained by the function $\boldsymbol{P}_{k}$ defined in \eqref{EQ4*}.
	\end{thm}
	\bpf
	Let  $f = h + \overline{g} \in \mathcal{S}_{H_q}^{0}(\mathcal{S},K)$ have the series representation given by \eqref{equation1***}. Then there exists some $\epsilon$ with $\vert \epsilon \vert = 1$ and a function $\psi \in \mathcal{S}$ such that
	$
	\phi_{\epsilon} = h + \epsilon g \prec_q \psi
	$. This implies there exist analytic functions $\xi(z)$ and $\eta(z)$  in the unit disk $\mathbb{D}$ satisfying $\vert \xi(z) \vert \leq  1$, $\vert \eta(z) \vert < 1$ for all $z \in \mathbb{D}$, with $\eta(0) = 0$ such that
	$$
	\phi_{\epsilon}(z) = \xi(z)\psi(\eta(z)) =z+ \sum_{n=2}^\infty {\epsilon}_nz^n.
	$$
	By Theorem \ref{thmA}, we have $\vert\epsilon_n\vert \leq n$, that is, $$\big \vert \vert a_n \vert - \vert b_n \vert \big\vert\leq \vert a_n+\epsilon b_n \vert \leq n ~\text{for all}~ n \geq 2.$$ Since the dilatation $\omega(z) = g'(z)/h'(z)$ satisfies $\omega(0) = 0$ and
	$\vert \omega(z) \vert \leq k$ for all $z \in \mathbb{D}$, it follows from
	Schwarz's lemma that
	$\vert \omega(z) \vert \leq k \vert z \vert$ for every $z \in \mathbb{D}$, and
	$\vert \omega'(0) \vert \leq k$.  This implies
	$$
	\frac{\omega(z)}{1 + \epsilon \,\omega(z)}
	\prec
	\frac{k z}{1 + \epsilon \,k  z}.
	$$
	As $\phi_{\epsilon} = h + \epsilon g$, we have $$h'(z)(1 + \epsilon\omega(z))=h'(z)+\epsilon g'(z)=\phi_{\epsilon}'(z)$$
	so that
	$$
	g'(z) = \omega(z) h'(z) = \dfrac{\omega(z)}{1 + \epsilon\omega(z)} \, \phi_{\epsilon}'(z).
	$$
	The Taylor coefficients of $ w(z)/(1 + \epsilon \omega(z))$  and   $\phi_\epsilon'(z)$ are dominated in modulus by those of $kz/(1-kz)$ and the derivative of the Koebe function $
	\boldsymbol{k}(z) = z/(1 - z)^2$, respectively. In accordance with the terminology from \cite[Chapter~7, p. 82]{Goodman:Univalent1}, we may express this as
	$$
	g^{\prime}(z)=\dfrac{\omega(z)}{1 + \epsilon \omega(z)} \, \phi_{\epsilon}'(z) \ll \frac{kz}{1 - kz} \cdot \frac{1 + z}{(1 - z)^3}  =
	\sum_{n=1}^{\infty} \left( \sum_{m=1}^{n} k^{\,m} (n-m+1)^{2} \right) z^{\,n}
	.
	$$
	Thus, as $b_{1} = 0$ and
	$
	g'(z) = \sum_{n=2}^{\infty} n b_{n} z^{\,n-1},
	$
	the above observation yields that
	$$
	n \,\vert b_{n} \vert \leq \sum_{m=1}^{n-1} k^{\,m} (n-m)^2
	\quad \text{for } n \geq 2,
	$$
	and thus (cf.~\cite{LiPonnusamy2025}) we obtain
	\beqq \label{eq11+}
	\vert b_n \vert \leq \frac{1}{n} \sum_{m=1}^{n-1} k^{\,m} (n-m)^2  = n\sum_{m=1}^{n-1} k^{\,m} +\frac{1}{n} \sum_{m=1}^{n-1} m^2 k^{\,m}-2\sum_{m=1}^{n-1} m k^{\,m} \quad \text{for } n \geq 2.
	\eeqq
	\noindent so that
	$
	\vert b_{n} \vert \leq  \boldsymbol{B}(n,k) .
	$
	Since $h(z) = \phi_\epsilon(z) - \epsilon g(z)$, it follows that
	$$
	\vert a_n \vert  \leq \big \vert \vert a_n \vert - \vert b_n \vert \big\vert + \vert b_n \vert \leq n + \boldsymbol{B}(n,k),
	$$
	which gives as in \cite{LiPonnusamy2025}, $\vert a_n \vert\le \boldsymbol{A}(n,k)$.  Next, we have
	$$
	h(re^{i\theta}) = \int_0^r h'(s e^{i\theta}) e^{i\theta} \, ds
	\quad \text{and} \quad
	g(re^{i\theta}) = \int_0^r g'(s e^{i\theta}) e^{i\theta} \, ds,
	$$
	and
	$$
	\vert f(z)\vert = \vert h(z) + \overline{g(z)}\vert \le \vert h(z)\vert + \vert g(z)\vert
	$$
	which imply
	\beq \label{eqn17+}
	\vert f(z)\vert \le \int_0^r \vert f_z(s e^{i\theta})\vert \, ds
	+ \int_0^r \vert f_{\bar{z}}(s e^{i\theta})\vert \, ds.
	\eeq
	Since $f = h + \overline{g} \in \mathcal{S}_{H_q}^{0}(\mathcal{S},K)$, it follows from Theorem~\ref{thmA} that
	$$
	\phi_{\epsilon}(z) = h(z) + \epsilon g(z) \ll \frac{z}{(1-z)^{2}}
	\qquad \text{for some } \epsilon \text{ with } \vert \epsilon \vert = 1.
	$$
	Furthermore, by \cite[p.~82, Theorem~5]{Goodman:Univalent1}, we have
	$$
	\phi_{\epsilon}'(z)\ll \frac{1+z}{(1-z)^{3}}
	\qquad \text{and} \qquad
	\vert \phi_{\epsilon}'(z)\vert \le \frac{1+r}{(1-r)^{3}},
	\quad \vert z \vert = r < 1.
	$$
	Next, as $f_z=h'$, $f_{\overline{z}}=\overline{g'}$ and $\vert w(z)\vert \le k\vert z \vert$, we obtain
	\beq \label{eqn18+}
	\vert f_z(z) \vert
	= \left\vert \frac{\phi_{\epsilon}'(z)}{1 +\epsilon \omega(z)} \right\vert
	\le \frac{\vert\phi_{\epsilon}'(z) \vert}{1 - k\vert  z \vert}\le \frac{1+r}{(1-rk)(1-r)^3},
	\eeq
	and
	\beq \label{eqn19+}
	\vert f_{\bar{z}}(z) \vert
	= \left\vert \frac{\omega(z)\phi_{\epsilon}'(z)}{1 +\epsilon \omega(z)} \right\vert
	\le \frac{kr(1+r)}{(1-rk)(1-r)^3}.
	\eeq
	Combining inequalities \eqref{eqn18+} and \eqref{eqn19+} with \eqref{eqn17+}, we obtain the desired inequality
	\beqq
	\vert f(z)\vert &\le&
	\int_0^r \frac{1+s}{(1-sk)(1-s)^3} \, ds
	+ \int_0^r \frac{ks(1+s)}{(1-sk)(1-s)^3} \, ds
	\\
	&=&
	\frac{r (k + 1)}{( 1-k)( 1-r)^2}- \frac{4 k r}{(1-k)^2 (1-r)}
	+ \frac{2 k (1 + k)}{(k - 1)^3}
	\log\left( \frac{ 1-r}{  1-k r} \right).
	\eeqq
	To show that the bounds are optimal, it suffices to verify that
	$\boldsymbol{P}_{k} \in \mathcal{S}_{H_q}^{0}(\mathcal{S},K)$.
	This follows directly from the construction of
	$\boldsymbol{P}_{k} $
	given in \cite{LiPonnusamy2025, Wang2024}, together with the fact that every function
	in the class $\mathcal{S}$ is subordinate to itself.
	This completes the proof.
	\epf
	\begin{rem}
		For the case $n=2$, Theorem~\ref{thm1+} gives the following sharp coefficient estimate for the class $\mathcal{S}_{H_q}^{0}(\mathcal{S},K)$:
		$$
		\vert a_2 \vert \le \frac{5K+3}{2K+2},
		$$
		which confirms Conjecture~1 of \cite{Wang2024} for this class.
		
	\end{rem}
	\begin{rem}
		As pointed out in \cite{LiPonnusamy2025},
		$\boldsymbol{A}(n,k)$ and $\boldsymbol{B}(n,k)$ are increasing functions of $k\in[0,1)$.
		If $k\to 1^{-}$, then Theorem \ref{thm1+} shows that
		$$
		\lim_{k\to 1^{-}} \boldsymbol{A}(n,k)=\frac{(2n+1)(n+1)}{6}
		\quad \text{and} \quad
		\lim_{k\to 1^{-}} \boldsymbol{ B}(n,k)=\frac{(2n-1)(n-1)}{6}.
		$$
		In this limiting case, the bounds obtained from Theorem~\ref{thm1+} coincide with those discussed in Conjecture~\ref{conj:A*} and Theorem~1 in \cite{Jasbir}.
	\end{rem}
	\begin{rem}
		For $0<r<1$,  we note that 
		$$
		\lim_{k \to 1^-} S(k,r)
		= \frac{r(r^2+3)}{3(1-r)^3}.
		$$
		Thus, the upper bound for $\lvert f(z)\rvert$ in Theorem~\ref{thm1+} coincides with the  corresponding upper bound for functions $f \in \mathcal{S}_{H_q}^0(\mathcal{S})$ as $k \to 1^{-}$ (see \cite[Theorem~1]{Jasbir}).
	\end{rem}

	
	We now state and prove our next result which determines the range of $p$ such that every $f \in \mathcal{S}_{H_q}(\mathcal{S}, K)$ belongs to 
	the Hardy space $\bf{h}^p$ and belong to the Bergman space $\mathbf{a}^{\mathbf{p}}_{\boldsymbol{\beta}}$, respectively, thereby  addressing Problem \ref{prob-hardy}.
	
	\begin{thm}\label{thm2@}
		Assume that $ f  \in  \mathcal{S}_{H_q}(\mathcal{S}, K)$.
		Then $ f\in \bf{h}^p$ for $0 < p < 1/2$, and $ f\in \mathbf{a}^{\mathbf{p}}_{\boldsymbol{\beta}}$ for $0<p<(2+\beta)/2$ and $\beta >-1$.
	\end{thm}
	\bpf
	Suppose that $f = h + \overline{g} \in \mathcal{S}_{H_q}(\mathcal{S},K)$.
	Then $c = g'(0)$  is a complex constant with $\vert c \vert < 1$ and that the function
	\beqq
	f_0 = \frac{f - \overline{c}\,\overline{f}}{1 - \vert c \vert^2}
	\in \mathcal{S}_{H_q}^{0}(\mathcal{S}, K_0),
	\eeqq
	such that $f_0$ is $K_0$-quasiconformal for some $K_0 \ge 1$ (see \cite[ p.~15]{LiPonnusamy2025}).
	We write $f_0 = h_0 + \overline{g_0}$ and set
	\beqq
	k_0 := \frac{K_0 - 1}{K_0 + 1}.
	\eeqq
	For all $0<r<1$, we have
	$$
	0<\frac{1-r}{1-k_0 r}<1
	\quad \Longrightarrow \quad
	\log\!\left(\frac{1-r}{1-k_0 r}\right)\le 0.
	$$
	Moreover, using the inequality $-\log (1-y) \le y/(1-y)$ for $y \in [0,1)$, we have
	$$
	-\log x \le \frac{1-x}{x}
	\quad \text{for}\quad 0 \le x <1
	$$
	and thus,
	$$
	\left\vert
	\log\!\left(\frac{1-r}{1-k_0 r}\right)
	\right\vert
	\le \frac{(1-k_0)r}{1-r}.
	$$
	Hence,
	$$
	\left\vert
	\frac{2k_0(1+k_0)}{(k_0-1)^3}
	\log\!\left(\frac{1-r}{1-k_0 r}\right)
	\right\vert
	\le
	\frac{C_1(k_0)\, r}{1-r},
	$$
	for a constant $C_1(k_0)>0$.
	Thus, there exists a constant $C=C(k_0)>0$ such that for all $0<r<1$,
	$$
	M(r)
	=
	\max_{0 \le \theta \le 2\pi}
	\vert f_0(re^{i\theta}) \vert\le
	\frac{C\,r}{(1-r)^2}~\qquad (\text{using Theorem~\ref{thm1+}}).
	$$
	Applying Lemma~\ref{lemB++}, we find that
	\beqq
	\lim_{r \to 1^-} M_p^p(r, f_0)
	&\le&
	\frac{2(1+k_0^2)\bigl(\vert p-2 \vert + 1\bigr)}{1-k_0^2}
	\int_0^1 M(s)^p s^{-1}\, ds\\
	&\le&
	N(p,k_0)
	\int_0^1 \frac{s^{p-1}}{(1-s)^{2p}}\, ds= N(p,k_0) B(p,1-2p),
	\eeqq
	where  $N(p,k_0)> 0$ is a constant depending on $p$ and $k_0$. Therefore, $f_0 \in {\bf{h}^p} \text{ for } 0<p < 1/2.$ It follows that $f - \overline{c}\,\overline{f} \in \bf{h}^p$. Using the inequality
	$$
	\| f - \overline{c}\,\overline{f} \|_p^p
	\ge \| f \|_p^p - \vert c \vert^p \| f \|_p^p
	= \bigl(1 - \vert c \vert^p\bigr)\, \| f \|_p^p,
	$$
	we conclude that $f \in \bf{h}^p$ for $0 < p < 1/2$.
	
	Next, following the proof of Theorem~3, we obtain
	\beqq
	\vert\vert f_0 \vert\vert_{\mathbf{a}^{\mathbf{p}}_{\boldsymbol{\beta}}}^p
	&\lesssim&
	\int_0^1 M(\rho)^p \rho^{-1}(1-\rho^2)^{\beta+1}\,d\rho\\
	&\lesssim&\int_{0}^{1}
	\frac{\rho^{p-1}}
	{(1-\rho)^{2 p}}(1-\rho^2)^{\beta+1}
	\, d\rho\\
	&\lesssim& \int_0^1 (1-\rho)^{-2 p+\beta+1}\, \rho^{p-1}\,d\rho = B(p, 2+\beta-2 p),
	\eeqq
	which shows that
	$f_0 \in \mathbf{a}^{\mathbf{p}}_{\boldsymbol{\beta}}$ whenever $0 < p < (\beta+2)/2$.
	\epf

	We now state and prove results for odd $K$-quasiconformal harmonic mappings discussed in Section \ref{section1.7}.
	\begin{thm}\label{thm5}
		Let $f\in \mathcal{S}_{H}^0(\mathcal{S},K)$ and that
		$ f(z)=z+\sum_{n=1}^{\infty} a_n z^{2n+1} +\overline{ \sum_{n=1}^{\infty} b_n z^{2n+1}}
		$ be an odd function. Then for all $n \ge 1$, we have the following inequalities:
		$$
		\big \vert \vert a_n \vert  - \vert b_n \vert \big \vert  < \lambda ,$$
		$$
		\vert a_n \vert  <
		\frac{\lambda \left[
			(2n+1) - (2n+3)k + k^{\,n+1} + k^{\,n+2}\right]}
		{(2n+1)(1-k)^{2}}
		,$$
		and
		$$
		\vert b_n \vert  <
		\frac{
			\lambda \,k\!\left[(2n-1) - (2n+1)k + k^{\,n} + k^{\,n+1}\right]
		}{
			(2n+1)(1-k)^{2}
		},
		$$
		where $\lambda = 1.1305$.
		Furthermore, $$\vert f(z)\vert \leq  T(k,r) ~\quad \text{for} ~\quad \vert z \vert = r < 1 ,$$ where
		\beqq
		T(k,r) := \frac{(k+1)r}{(1-k)(1-r^{2})}
		- \frac{2k}{(1-k)^{2}} \log \frac{1+r}{1-r}
		+ \frac{\sqrt{k}(k+1)}{(1-k)^{2}} \log\!\left(\frac{1+\sqrt{k}\,r}{1-\sqrt{k}\,r}\right).
		\eeqq
		Moreover, the bound for $\vert f(z)\vert$ is sharp for the function $F_{k}(z)=H_k(z)+\overline{G_k(z)}$, where
		$$
		H_k(z)
		= \frac{z}{(1-k)(1-z^{2})}
		\;-\; \frac{k}{(1-k)^{2}} \log\!\frac{1+z}{1-z}
		\;+\; \frac{\sqrt{k}(k+1)}{2(1-k)^{2}}
		\log\!\frac{1+\sqrt{k}\,z}{1-\sqrt{k}\,z},
		$$
		and
		$$
		G_k(z)
		= \,\frac{k\,z}{(1-k)(1-z^{2})}
		\;-\; \frac{k}{(1-k)^{2}} \log\!\frac{1+z}{1-z}
		\;+\; \frac{\sqrt{k}(k+1)}{2(1-k)^{2}}
		\log\!\frac{1+\sqrt{k}\,z}{1-\sqrt{k}\,z}.
		$$
	\end{thm}
	\bpf
	Let $f\in \mathcal{S}_{H}^{0}(\mathcal{S},K)$ and
	$$ f(z)=h(z)+ \overline{g(z)}=z+\sum_{n=1}^{\infty} a_n z^{2n+1} +\overline{\sum_{n=1}^{\infty} b_n z^{2n+1}}  .
	$$
	Then $\phi_{\epsilon} = h + \epsilon g \in \mathcal{S}$ for some $\epsilon$ such that $\vert \epsilon \vert = 1$,  and
	$$
	\phi_{\epsilon}(z)  = z + \sum_{n=2}^{\infty} \epsilon_{n} z^{2n+1} ~\mbox{ with }~ \epsilon_{n}=a_n+\epsilon b_n .
	$$
	By the  Theorem \ref{thmA1}, $\vert \epsilon_{n} \vert < \lambda$ for all $n \ge 1$ and hence, $\big \vert \vert a_n \vert - \vert b_n \vert \big\vert\leq \vert a_n+\epsilon b_n \vert < \lambda$ for all $n \geq 1$. It follows that the dilatation $\omega(z) = g'(z)/h'(z)$ is an even function  such that $\omega(0) = 0=w'(0)$ and
	$\vert \omega(z) \vert \leq k$ for all $z \in \mathbb{D}$. Therefore,
	$\vert \omega(z) \vert \leq k \vert z \vert^2$ for $z \in \mathbb{D}$.
	Let
	$$
	\frac{\omega(z)}{1+\epsilon\,\omega(z)}
	= \sum_{n=1}^{\infty} \omega_n z^{2n}.$$
	Since $\vert \omega(z) \vert \leq k \vert z \vert^2$, we have
	$$
	\frac{ \omega(z)}{1+\epsilon\,\omega(z)}  \ll \frac{kz^2}{1-kz^2}=\sum_{n=1}^{\infty} k^{n} z^{2n},
	$$
	from which we obtain $
	\vert \omega_n \vert \le k^n ~~ \text{for all } n\ge 1.$ Next, we can express $g'(z)$ as follows:
	\beqq
	g'(z)
	= \frac{\phi'_{\epsilon}(z)\,\omega(z)}{1+\epsilon\,\omega(z)}
	&=& \left( \epsilon_0 + \sum_{n=1}^{\infty} (2n+1)\epsilon_n z^{2n} \right)
	\left( \sum_{n=1}^{\infty} \omega_n z^{2n} \right)\\
	&=&
	\sum_{n=1}^{\infty}
	\left(
	\epsilon_0 \,\omega_{n}
	+
	\sum_{l=1}^{n-1} (2k+1)\,\epsilon_{l}\,\omega_{n-l}
	\right)
	z^{2n},
	\eeqq
	where $\epsilon_0 = 1$. From this, we conclude that for all $n \ge 1$,
	\beqq
	(2n+1)\,\vert b_n \vert
	&\le& \vert \omega_n \vert
	+ \sum_{l=1}^{n-1} (2l+1)\,\vert \epsilon_l \vert\, \vert \omega_{n-l}\vert  \\
	&<& \vert \omega_n \vert
	+ \lambda \sum_{l=1}^{n-1} (2l+1)\, \vert \omega_{n-l} \vert \qquad (\text{since}~~ \vert \epsilon_l \vert < \lambda~~ \text{for all}~~l \ge 1)\\
	&=& \vert \omega_n \vert
	+ \lambda \sum_{m=1}^{n-1} (2(n-m)+1)\, \vert \omega_{m} \vert\\
	&=&  (1- \lambda)\vert \omega_n \vert+ \lambda \sum_{m=1}^{n} (2(n-m)+1)\, \vert \omega_{m} \vert\\
	&\le& \lambda \sum_{m=1}^{n} k^m(2(n-m)+1)\,\qquad (\text{since} ~~\vert \omega_n \vert \le k^n ~~ \text{for all } n\ge 1).
	\eeqq
	That is,
	$$
	\vert b_n \vert  <
	\frac{
		\lambda \,k\!\left[(2n-1) - (2n+1)k + k^{\,n} + k^{\,n+1}\right]
	}{
		(2n+1)(1-k)^{2}
	}.
	$$
	From the equation \(h(z)=\phi_{\epsilon}(z)-\epsilon\, g(z)\), we obtain
	$$
	\vert a_n \vert
	\le \vert \epsilon_n \vert + \vert b_n \vert
	< \frac{\lambda \left[
		(2n+1) - (2n+3)k + k^{\,n+1} + k^{\,n+2}\right]}
	{(2n+1)(1-k)^{2}},
	$$
	for all  $n \ge 1$.
	By invoking Lemma \ref{Lemma5A}, we obtain
	\beq \label{eqn17**}
	\vert f_z(z) \vert
	= \left\vert \frac{\phi_{\epsilon}'(z)}{1 +\epsilon \omega(z)} \right\vert
	\le \frac{\vert\phi_{\epsilon}'(z) \vert}{1 - k\vert  z \vert^2}\le \frac{1+r^2}{(1-kr^2)(1-r^2)^2},
	\eeq
	and thus,
	\begin{equation}\label{eqn18**}
		\vert f_{\bar z}(z) \vert
		= \left\vert \frac{\omega(z)\,\phi_\epsilon'(z)}{1+\epsilon\,\omega(z)} \right\vert
		\le  \frac{k r^2(1+r^2)}{(1-k r^2)(1-r^2)^2}.
	\end{equation}
	
	The inequalities \eqref{eqn17**} and \eqref{eqn18**}  with \eqref{eqn17+} can now be used to obtain the desired inequality, namely,
	\beqq
	\vert f(z)\vert &\le&
	\int_0^r \frac{1+s^2}{(1-ks^2)(1-s^2)^2} \, ds
	+ \int_0^r \frac{ks^2(1+s^2)}{(1-ks^2)(1-s^2)^2} \, ds
	\\
	&=&
	\frac{(k+1)r}{(1-k)(1-r^{2})}
	- \frac{2k}{(1-k)^{2}} \log \frac{1+r}{1-r}
	+ \frac{\sqrt{k}(k+1)}{(1-k)^{2}} \log\!\left(\frac{1+\sqrt{k}\,r}{1-\sqrt{k}\,r}\right).
	\eeqq
	Next, we examine the sharpness of the obtained bound. Consider
	$ h+ g=\phi$,  where $\phi(z)=z/(1 - z^{2})$  is the square-root transform of the
	Koebe function $\boldsymbol{k}(z) = z/(1 - z)^2$, and  $\phi$ is clearly univalent and convex in the vertical
	direction in $\mathbb{D}$ and therefore, close-to-convex in $\mathbb{D}$. By Theorem~\ref{clunie}, the harmonic mapping
	$f = h + \overline{g}$ is univalent and convex in the vertical direction,  provided $f$ is locally univalent in $\mathbb{D}$.
	In fact, the local univalency of $f$ can be assured by prescribing the dilatation
	$\omega(z) = -kz^{2}$, which yields the relation
	$
	g'(z) = -kz^{2} h'(z)
	$
	in $\mathbb{D}$. Differentiating
	$h + g =  z/(1-z^2)$ together with the relation $
	g'(z) = -kz^{2} h'(z)
	$ leads to the following pair of differential equations:
	$$
	h'(z) + g'(z) = \frac{1+z^{2}}{(1-z^{2})^{2}}
	\quad \text{and} \quad
	kz^{2} h'(z) + g'(z) = 0.
	$$
	Solving these yields
	$$
	h'(z) = \frac{1 + z^{2}}{(1-kz^2)(1 - z^{2})^{2}}
	\quad \text{and} \quad
	g'(z) = \frac{-kz^{2}(1 + z^{2})}{(1-kz^2)(1 - z^{2})^{2}}.
	$$
	Integration gives
	\beqq
	h(z)
	&=& \frac{z}{(1-k)(1-z^{2})}
	\;-\; \frac{k}{(1-k)^{2}} \log\!\frac{1+z}{1-z}
	\;+\; \frac{\sqrt{k}(k+1)}{2(1-k)^{2}}
	\log\!\frac{1+\sqrt{k}\,z}{1-\sqrt{k}\,z}\\
	&=& z+ \sum_{n=1}^{\infty}
	\left[
	\frac{1}{1-k}
	- \frac{2k}{(1-k)^2(2n+1)}
	+ \frac{k(k+1)}{(1-k)^2(2n+1)}\,k^{n}
	\right] z^{2n+1}\\
	&=& z+ \sum_{n=1}^{\infty}
	\left[
	\frac{(2n+1)-(2n+3)k+k^{n+1}+k^{n+2}}{(2n+1)(1-k)^{2}}
	\right] z^{2n+1}
	\eeqq
	and
	\beqq
	g(z)
	&=& \,\frac{k\,z}{(1-k)(1-z^{2})}
	\;-\; \frac{k}{(1-k)^{2}} \log\!\frac{1+z}{1-z}
	\;+\; \frac{\sqrt{k}(k+1)}{2(1-k)^{2}}
	\log\!\frac{1+\sqrt{k}\,z}{1-\sqrt{k}\,z}\\
	&=& \sum_{n=1}^{\infty}
	\left[
	\frac{k}{1-k}
	- \frac{2k}{(1-k)^2(2n+1)}
	+ \frac{k(k+1)}{(1-k)^2(2n+1)}\,k^{n}
	\right] z^{2n+1}\\
	&=& \sum_{n=1}^{\infty}
	\left[
	\frac{k\big[(2n-1)-(2n+1)k+k^{n}+k^{\,n+1}\big]}{(2n+1)(1-k)^{2}}
	\right] z^{2n+1}.
	\eeqq
	The construction and Theorem~\ref{clunie} show that
	$f = h + \overline{g}$ belongs to the class $\mathcal{S}_{H}^{0}(\mathcal{S},K)$.
	The proof is complete.
	\epf
	\begin{rem}\label{rem4}
		Let $f \in  \mathcal{S}_{H}^0(\mathcal{S})$ be an odd function.  Then allowing $k\to 1^{-}$ in Theorem \ref{thm5} gives
		$$
		\vert f(z) \vert \le \lim_{k \to 1^{-}} T(k,r)
		= \frac{1}{2}
		\left[\frac{1}{2}\log \frac{1+r}{1-r}
		+ \frac{r(1+r^{2})}{(1-r^{2})^{2}} \right]~\quad \text{for} ~\quad \vert z \vert = r < 1
		$$\
		Moreover, the bound is sharp for the function (cf. \cite[Example 1]{Jaglan})
		\beqq \label{8-}
		f(z)
		= \frac{1}{8}\log\!\left(\frac{1+z}{1-z}\right)
		- \frac{z^{3}-3z}{4(1-z^{2})^{2}}
		+\overline{ \frac{1}{8}\log\!\left(\frac{1-z}{1+z}\right)
			- \frac{3z^{3}-z}{4(1-z^{2})^{2}}}.
		\eeqq
		Note that the growth estimate for odd functions $f \in \mathcal{S}_H^{0}(\mathcal{S}, K) $ coincides with
		the growth estimate for odd functions $f \in \mathcal{S}_{H}^0(\mathcal{S})$ as $k \to 1^{-}$
	\end{rem}
	
	\begin{thm}\label{thm8++}
		Let $ f \in \mathcal{S}_{H}(\mathcal{S},K) $ be an odd function.
		Then $f \in \bf{h}^{p}$ for $0 < p < 1$ and  $f \in \mathbf{a}^{\mathbf{p}}_{\boldsymbol{\beta}}$ for $0<p<(2+\beta)$.
	\end{thm}
	\bpf
	The proof follows from the methodology similar to that of Theorem~\ref{thm2@}. For completeness, we briefly outline some of the main steps.
	\noindent Recall the  inequality,
	\beq\label{20--}
	\log \frac{1+x}{1-x} \le \frac{2x}{1-x^2}, ~\quad \text{for} \quad x \in [0,1),
	\eeq
	which follows by a comparison of the coefficients of the corresponding power series on both sides.
	Setting $x=r$ and $x=\sqrt{k_0}\,r$ in \eqref{20--} yields respectively the inequalities
	$$
	\log \frac{1+r}{1-r} \le \frac{2r}{1-r^2}\le \frac{2r}{1-r},
	$$
	and
	$$
	\log\!\left(\frac{1+\sqrt{k_0}\,r}{1-\sqrt{k_0}\,r}\right)
	\le \frac{2\sqrt{k_0}\,r}{1-k_0 r^2}
	\le \frac{2\sqrt{k_0}\,r}{1-r^2}\le \frac{2\sqrt{k_0}\,r}{1-r},
	$$
	since $k_0 \in (0,1)$ and $r\in (0,1)$. Thus, there exists a constant $C=C(k_0)>0$ such that for all $0<r<1$
	\beqq
	M(r)
	=
	\max_{0 \le \theta \le 2\pi}
	\vert f(re^{i\theta}) \vert \leq \,\frac{Cr}{1-r}~\qquad (\text{using Theorem~\ref{thm5}}).
	\eeqq
	Now, invoking Lemma~\ref{lemB++}, we obtain the desired result.
	\epf

	%
	%

	\subsection*{Acknowledgement}
	This work was partially supported by the FIST program of the Department of Science and Technology, Government of India (Reference No.~SR/FST/MS-I/2018/22(C)). The first author acknowledges the Council of Scientific and Industrial Research (CSIR), India, for financial support through a Senior Research Fellowship. The work of the second author is supported by Mathematical Research Impact Centric Support (MATRICS) grant (File No.: ANRF/ARGM/2025/\\000370/MTR) by the Anusandhan National Research Foundation (ANRF), Government of India. The third author was partially supported by the Core Research Grant (CRG/2022/008920) from the Anusandhan National Research Foundation (ANRF), Government of India.

	\subsection*{Conflict of Interest Statement}
	The authors declare that they have no conflict of interest, regarding the publication of this paper.
	
	\subsection*{Data Availability Statement}
	The authors declare that this research is purely theoretical and does not associate with any datas.
	
\end{document}